\documentstyle[11pt,a4wide]{article}

\font \tenmsb=msbm10 scaled \magstep 1
\font \sevenmsb=msbm7 scaled \magstep 1
\font \fivemsb=msbm5 scaled \magstep 1
\newfam \msbfam
\textfont \msbfam=\tenmsb
\scriptfont \msbfam=\sevenmsb
\scriptscriptfont \msbfam=\fivemsb
\def \Bbb#1{\fam \msbfam \relax#1}

\font \teneufm=eufm10 scaled \magstep 1
\font \seveneufm=eufm7 scaled \magstep 1
\font \fiveeufm=eufm5 scaled \magstep 1
\newfam \eufmfam
\textfont \eufmfam=\teneufm
\scriptfont \eufmfam=\seveneufm
\scriptscriptfont \eufmfam=\fiveeufm
\def \frak#1{{\fam \eufmfam \relax#1}}

\title{\bf A q-ANALOGUE OF THE BEREZIN QUANTIZATION OF THE UNIT DISC}

\author{\sl D. Shklyarov$^\diamond$ \and \sl S.  Sinel'shchikov$^\star$
\thanks{Partially supported by the grant INTAS-RFBR-95-418} \and \sl L.
Vaksman$^\star$ \thanks{Partially supported by the grant INTAS-94-4720}}

\date{\tt $^\star$Institute for Low Temperature Physics \& Engineering\\
National Academy of Sciences of Ukraine\\ \ \\
$^\diamond$Kharkov State University\\
Department of Mechanics and Mathematics}

\begin{document}

\maketitle

 It is well known \cite{MO-N, CGR} that the Berezin quantization can be used
to construct a formal deformation of the algebra $C^\infty(U)$ of smooth
functions in the unit disc $U \subset{\Bbb C}$.

 In the well known work \cite{KL} a family of q-analogues of the unit disc $U
\subset{\Bbb C}$ was introduced. The simplest, and hence the most important
one among those was also considered in \cite{NN, F, CK}. With these results
as a background, we produce an explicit formula for a formal deformation of
the quantum disc. To develop it, a Berezin method \cite{B} is used.

 Consider an involutive algebra ${\rm Pol}({\Bbb C})_q$, determined by its
generator $z$ and the relation
$$z^*z=q^2zz^*+1-q^2.\eqno(1)$$
This $*$-algebra may be treated as a q-analogue of the Weyl algebra
\cite{CK}. However, we use an approach of \cite{KL, NN, F} to the above
algebra as a q-analogue of the polynomial algebra on the plane. We use a
standard differential calculus \cite{Zu}. Specifically, let $\Omega_q({\Bbb
C})$ be an involutive algebra given by its generators $z$, $dz$, and the
relations (1),
$$dz \cdot z^*=q^{-2}z^* \cdot dz,\qquad dz \cdot z=q^2z \cdot dz,\qquad dz
\cdot dz^*=-q^{-2}dz^* \cdot dz,\qquad dz \cdot dz=0.$$

 A differential $d:\Omega_q({\Bbb C})\to \Omega_q({\Bbb C})$ is defined in a
standard manner on the elements $z$, $z^*$, $dz$, $dz^*$, and then extended
in a unique way up to a differentiation of the graded algebra
$\Omega_q({\Bbb C})$ (${\rm deg}(z)={\rm deg}(z^*)=0$, ${\rm deg}(dz)={\rm
deg}(dz^*)=1$). Partial derivatives ${\textstyle \partial^{(l)}\over
\textstyle \partial z}$, ${\textstyle \partial^{(r)}\over \textstyle
\partial z}$, ${\textstyle \partial^{(l)}\over \textstyle \partial z^*}$,
${\textstyle \partial^{(r)}\over \textstyle \partial z^*}$ are defined as
such linear operators in ${\rm Pol}({\Bbb C})_q$ that for all $f \in{\rm
Pol}({\Bbb C})_q$
$$df=dz \cdot{\partial^{(l)}f \over \partial z}+dz^*{\partial^{(l)}f \over
\partial z^*}={\partial^{(r)}f \over \partial z}\cdot dz+{\partial^{(r)}f
\over \partial z^*}\cdot dz^*.$$
The bilinear maps $L:{\rm Pol}({\Bbb C})_q \times{\rm Pol}({\Bbb C})_q
\to{\rm Pol}({\Bbb C})_q$ of the form
$$L:f_1 \times f_2 \mapsto \sum_{i,j,k,m=0}^{N(L)}c_{ijkm}\cdot
\left(\left({\partial^{(r)}f \over \partial z^*}\right)^if_1
\right)\cdot z^{*j}z^k \cdot \left(\left({\partial^{(l)}f \over \partial
z}\right)^mf_2 \right),$$
with $c_{ijkm}\in{\Bbb C}$, will be called q-bidifferential operators.

 Our goal is to produce a formal deformation of the multiplication law in
${\rm Pol}({\Bbb C})_q$
$$*:{\rm Pol}({\Bbb C})_q \times{\rm Pol}({\Bbb C})_q \to{\rm Pol}({\Bbb
C})_q[[t]],$$
$$*:f_1 \times f_2 \mapsto f_1 \cdot f_2+\sum_{j \in{\Bbb Z}_+}t^j \cdot
C_j(f_1,f_2),$$
and then to get explicit formulae for q-bidifferential operators $\{C_j\}_{j
\in{\Bbb Z}_+}$. (The requirement of formal associativity \cite{BFFLS} will
be satisfied.)

 It should be noted that the operators ${\textstyle \partial^{(l)}\over
\textstyle \partial z}$, ${\textstyle \partial^{(l)}\over \textstyle
\partial z^*}$ ${\textstyle \partial^{(r)}\over \textstyle \partial z}$,
${\textstyle \partial^{(r)}\over \textstyle \partial z^*}$ are extended "by
a continuity" from the space of polynomials ${\rm Pol}({\Bbb C})_q$ onto the
space $D(U)^\prime$ of series of the form $f=\displaystyle \sum \limits_{j,k
\in{\Bbb Z}_+}a_{jk}z^jz^{*k}$, $a_{jk}\in{\Bbb C}$. This allows one to
extend the $*$-product described in the sequel from the polynomial algebra
onto the "algebra of smooth functions in the quantum disc".

 Everywhere below $q,t \in(0,1)$.

 We follow \cite{KL} in considering the operator of "weighted shift"
$\widehat{z}$ given by $\widehat{z}e_m=e_{m+1}$ with respect to an
orthogonal base such that
$$\|e_m \|^2=\frac{(q^2;q^2)_m}{(q^2t;q^2)_m}.\eqno(2)$$
A standard notation $(a;q^2)_m=(1-a)(1-q^2a)\cdot \ldots
\cdot(1-q^{2(m-1)}a)$ is used here.

 It was explained in \cite{KL} that the operators $\widehat{z}$,
$\widehat{z}^*$ supply a two-parameter quantization of the unit disc in the
complex plane.

 A covariant symbol of the operator $\widehat{f}=\displaystyle \sum
\limits_{i,j}a_{ij}\widehat{z}^i \widehat{z}^{*j}$, $a_{ij}\in{\Bbb C}$ is
defined to be the series $f=\displaystyle \sum \limits_{ij}a_{ij}z^iz^{*j}$.
The map $f \mapsto \widehat{f}$ is a q-analogue of the Berezin
quantization.\footnote{In this work we consider only polynomial symbols $f
\in{\rm Pol}({\Bbb C})_q$. The general case was investigated in
\cite{SSV5}.} The formal deformation $*$ of the multiplication in ${\rm
Pol}({\Bbb C})_q$ is defined via this quantization:
$$\widehat{f_1*f_2}=\widehat{f_1}\cdot \widehat{f_2},\qquad f_1,f_2 \in{\rm
Pol}({\Bbb C})_q.$$
 This deformation is well defined, as one can see from the following
observation. (2) implies the relation
$$\widehat{z}^*\widehat{z}=q^2 \widehat{z}\widehat{z}^*+1-q^2+t
\frac{1-q^2}{1-t}(1-\widehat{z}\widehat{z}^*)(1-\widehat{z}^*\widehat{z}).
\eqno(3)$$
Explicit formulae for the multiplication $*$ are also accessible via
multiple application of (3). However, in our opinion, the proof of those
formulae presented in \cite{SSV5} is more appropriate.

\medskip

Our principal result is

  {\it For all $f_1,f_2 \in{\rm Pol}({\Bbb C})_q$
$$f_1*f_2=(1-t)\sum_{j \in{\Bbb
Z}_+} t^jm(p_j(\stackrel{\sim}{\Box})f_1 \otimes f_2),$$
with \\
i) $m:{\rm Pol}({\Bbb C})_q \otimes {\rm Pol}({\Bbb C})_q \to{\rm Pol}({\Bbb
C})_q$ being a multiplication in ${\rm Pol}({\Bbb C})_q$,\\
ii) $\stackrel{\sim}{\Box}=q^{-2}(1-(1+q^{-2})z^*\otimes
z+q^{-2}z^{*2}\otimes z^2)\cdot{\textstyle \partial^{(r)}\over
\textstyle \partial z^*}\otimes{\textstyle \partial^{(l)}\over
\textstyle \partial z}$,\\
iii) $p_j(\stackrel{\sim}{\Box})= \displaystyle \sum
\limits_{k=0}^j
\frac{\textstyle(q^{-2j};q^2)_k}{\textstyle(q^2;q^2)^2_k}q^{2k} \cdot
\displaystyle \prod \limits_{i=0}^{k-1}\left(1-q^{2i}\left((1-q^2)^2
\stackrel{\sim}{\Box}+1+q^2 \right)+q^{4i+2}\right)$.}

\medskip

 The proof of this theorem is presented in \cite{SSV5}. It involves the
notions of a contravariant symbol and of a Berezin transform \cite{UU}. The
contravariant symbol of the operator $\widehat{f}=\displaystyle
\sum \limits_{i,j=1}^{N(\widehat{f})}a_{ij}\widehat{z}^{*i}\widehat{z}^j$,
$a_{ij}\in{\Bbb C}$, is defined to be a polynomial
$\stackrel{\circ}{f}=\displaystyle \sum
\limits_{i,j=1}^{N(\widehat{f})}a_{ij}z^{*i}z^j \in{\rm Pol}({\Bbb C})_q$.
The Berezin transformation $B_q:{\rm Pol}({\Bbb C})_q \to{\rm Pol}({\Bbb
C})_q[[t]]$, $B_q:\stackrel{\circ}{f}\mapsto f$ takes the contravariant
symbols of operators to their covariant symbols.

 Just as in the case of an ordinary disc, the proof of the theorem reduces
to constructing an "asymptotic expansion" of the Berezin transform. The
explicit form for the coefficients of this expansion could be found via
an application of the methods of the quantum group theory \cite{CP}.
($*$-algebra ${\rm Pol}({\Bbb C})_q$ is a $U_q \frak{sl}_{1,1}$-module
algebra, and the Berezin transformation is a morphism of $U_q
\frak{sl}_{1,1}$-modules.)

\end{document}